\documentclass{article}
\usepackage{graphicx}
\RequirePackage[OT1]{fontenc}
\RequirePackage{amsthm,amsmath,amssymb}
\RequirePackage[numbers]{natbib}
\RequirePackage[colorlinks,citecolor=blue,urlcolor=blue]{hyperref}
\usepackage[usenames]{color}

\numberwithin{equation}{section}
\theoremstyle{plain}

\newtheorem{assumption}{Assumption}[section]
\newtheorem{lemma}[assumption]{Lemma}
\newtheorem{definition}[assumption]{Definition}
\newtheorem{proposition}[assumption]{Proposition}

\newtheorem{theorem}[assumption]{Theorem}

\begin{document}

\title{Weak consistency of Markov chain Monte Carlo methods
}
\date{}
\author{Kengo KAMATANI\footnote{
Graduate School of Engineering Science, Osaka University
Machikaneyama-cho 1-3, Toyonaka-si, Osaka, 560-0043, Japan, kamatani@sigmath.es.osaka-u.ac.jp }
\footnote{Supported in part by Grant-in-Aid for JSPS Fellows (19-3140)
and  Grant-in-Aid for Young Scientists (B) 22740055.}
}

\maketitle

\begin{abstract}
%BEGIN update 20110213
Markov chain Monte Calro methods (MCMC) are commonly used in Bayesian statistics. 
In the last twenty years, many results have been established for the calculation of the exact
 convergence rate 
of MCMC methods. 
We introduce  another  rate of convergence for MCMC methods by approximation techniques. 
This rate can be obtained by the convergence of the Markov chain
to a diffusion process. 
We apply it to a  simple mixture model and obtain its convergence rate. 
Numerical simulations are performed to illustrate the effect of the rate. 
%END update 20110213
\end{abstract}

\section{Introduction}
Markov chain Monte Carlo (MCMC) method has become an essential tool
in any study that has a complicated posterior calculation problem.
Various new MCMC methods have been developed in the last two decades.
Theoretical support of this strategy has also been developed 
such as \citet{RT}, \citet{MT2} and many others.  
In particular, it was shown that the usual MCMC method produces an ergodic Markov chain (see \citet{TierneyAOS94} and \citet{RR}).

In practice, it is of great interest to study the convergence speed of these Markov chains. 
Various  quantitative bounds have been developed from the spectral approach by such as 
\citet{MR1097463} and  \citet{DKS}, and from the so-called (double) drift condition approach by 
\citet{MT1994} and \citet{R, ECP1054}. 
For an ergodic  Markov chain $\{X_m;m\ge 0\}$ on the state space $(E,\mathcal{E})$ with the  transition kernel $P$, they calculated the upper bound of 
\begin{equation}\label{tv}
\|\mathcal{L}(X_m)-\Pi\|=2\sup_{A\in \mathcal{E}}|\mathbb{P}(X_m\in A)-\Pi(A)|
\end{equation}
where $\Pi$ is the invariant distribution and $\mathcal{L}(X_m)$ is the law of $X_m$.  
In the former approach, if we can calculate the eigenvalues and the eigenfunctions for $P:L^2(\Pi)\rightarrow L^2(\Pi)$, then it is possible to calculate the almost exact bounds. On the other hand, although the latter approach does not provide tight bound, it is relatively easy to apply. 

To compare different MCMC methods, the above approaches may have difficulties, since we need to calculate  tight (upper and lower) bounds for two or more MCMC methods. However, without calculating such bounds, sometimes it is possible to compare different MCMC methods  by the asymptotic variance $\sigma^2_f$ in the following limit in $M\rightarrow\infty$:
\begin{equation}\nonumber
\sqrt{M}(I_M-I)\Rightarrow N(0,\sigma^2_f),\ I_M=\frac{1}{M}\sum_{m=0}^{M-1}f(X_m),\ I=\int_E f(x)\Pi(dx). 
\end{equation}
For this comparison, it is sufficient to show positivity of an operator in $L^2(\Pi)$ sense. 
This approach was studied in \citet{PESKUN01121973}, and later developed by \citet{zbMATH01319841} and \citet{mira1998ordering}. 
Although the application area of this is limited, this approach is particularly useful for the comparison of the so-called data-augmentation (DA) procedure
with its parameter-expanded extension (see \citet{HobertMarchev08}). 

These analysis  on MCMC procedures obtain the exact bound of the convergence rate or the exact comparison of MCMC procedures. 
We took a different approach in  \citet{Kamatani10}. 
Usually, MCMC procedures are complicated that prevent us from exact analysis. 
On the other hand, by using approximation theory, such as the  traditional large sample theory, 
sometimes it is easy to perform theoretical comparison among MCMC procedures. 
For this approximation, we introduce an index $n$, which tends to $\infty$. 
As $n\rightarrow\infty$,  if the following holds for any $M_n\rightarrow\infty$, the MCMC 
procedure is said to have the consistency in \cite{Kamatani10}:
\begin{equation}\label{lc}
I^n_{M_n}-I^n=o_{\mathbb{P}}(1),\ I^n_M=\frac{1}{M}\sum_{m=0}^{M-1}f(X_m^n),\ I^n=\int_E f(x)\Pi_n(dx)
\end{equation}
where $\{X_m^n;m\ge 0\}$ are sequences of $\Pi_n$-invariant  Markov chains generated by MCMC procedures. 
By Theorem 1 of \citet{Kamatani10}, under some regularity conditions, the DA procedure   satisfies this property. 
In practice,  if an MCMC procedure has the  consistency, it works fairly well. On the other hand, many popular MCMC methods do not satisfy this good convergence property
but satisfy a bad property
\begin{equation}\label{dc}
I^n_M-I^n_1=o_{\mathbb{P}_n}(1) 
\end{equation}
for any fixed $M\in\mathbb{N}$. This property means that the Monte Carlo estimation using $M$ iteration is no more helpful than that using only one iteration. 
Therefore we can classify MCMC procedure into two categories (\ref{lc}) or (\ref{dc}). 
Although these two categories do not cover all of the cases, this classification is useful in practice. 
However it does not tell the rate of convergence. 

In this paper, we introduce a further step of this approach. 
As mentioned earlier, the rate of convergence is useful to predict sufficient number of iteration until convergence, or 
to compare different MCMC procedures in details. 
We call $r_n$ the order of the weak consistency if (\ref{lc}) is satisfied for any $M_n\rightarrow\infty$
such that $M_n/r_n\rightarrow\infty$. If the MCMC procedure has the  consistency, 
we can take $r_n=1$. On the other hand, the order can be high if the performance of MCMC procedure is poor, that is, 
the condition (\ref{dc}) is satisfied. 
The order $r_n$ can be interpreted as the order of the sufficient number of iteration.

As an example we will consider the 
DA procedure for a simple mixture model $pN(\epsilon,1)+(1-p)N(0,1)$ for unknown $p$ but for known $\epsilon$. 
Since the performance of the DA procedure heavily depends on the parameter $\epsilon$, we let $\epsilon\rightarrow 0$ to illustrate the effect. 
This DA procedure works quite poorly if the true model is close to $N(0,1)$. 
The index $n$ is the sample size. 
It has the order $r_n=\epsilon^{-1}n^{1/2}$ and this shows the effects of both $\epsilon$ and the sample size $n$. 
This result comes from the fact that the trajectory of the DA procedure tends to a path of the 
stochastic process defined by
\begin{equation}\label{sde}
dX_t=(\alpha_1+X_t z-X_t^2 I)dt+\sqrt{2X_t}dW_t
\end{equation}
where  $I$ corresponds to the Fisher information matrix 
and $z$ corresponds to the scaled maximum likelihood estimator
(see Theorem \ref{gibbsweak}). 
It is probably well recognized that the trajectory of poor behaved MCMC procedure  looks like
a path of a diffusion process. This result is the first validation for this observation.

This paper is organized as follows. 
Section \ref{Sec2} we define (local) weak consistency. 
In Section \ref{GSM} we apply this to the simple mixture model. 
 Numerical results is provided in Section \ref{NR} which shows the 
effect of the order of the weak consistency. 

%\subsection{Notation}
%$\mathbf{R}^+=[0,\infty)$, $\mathbf{N}_0=\{0,1,\ldots\}$. 
%$\Theta$ is a Polish space
%having metric $d$ such that $d(a,b)\le 1$. 
%For measurable space $(X,\mathcal{X})$, $\mathcal{P}(X)$ 
%is the space of probability measures on $X$. 
%Probability transition kernel (Ptk) from $(X,\mathcal{X})$ to $(Y,\mathcal{Y})$ is 
%a function $K(dy;x)$ such that 
%$K(A;\cdot)$ is $\mathcal{X}$-measurable for any $A\in\mathcal{Y}$
%and $K(\cdot;x)\in\mathcal{P}(Y)$ for any $x\in X$. 
%We may write ``$P$ is a Ptk from $X$ to $Y$'' for short. 
%
%We write $\|\nu\|=\sup_{A\in\mathcal{X}}|\nu(A)|$
%for the total variation of a signed measure $\nu$. 
%Denote $\mathrm{BL}_1(\Theta)$ for a set of functions $\psi$ on $\Theta$ such that 
%$|\psi(a)-\psi(b)|\le d(a,b)$ and $|\psi(a)|\le 1$. 
%For two measure $\nu$ and $\mu$, 
%\[
%\sup_{\psi\in \mathrm{BL}_1(\Theta)}|\nu(\psi)-\mu(\psi)|
%\]
%is called a Bounded Lipshitz metric denoted by $w(\mu,\nu)$. 

\section{Local weak consistency of MCMC}\label{Sec2}
We write $[x]$ for the integer part of $x\in\mathbb{R}$. 

\subsection{Definition of local weak consistency}

In this section, we review  the (local) consistency and degeneracy
and also, we define the order of the weak consistency. 
Let $\Theta=\mathbb{R}^d$ be a parameter space. 
Suppose that observation $x_n$ is an element of a set $X_n$, and we are interested in the approximation of 
the posterior distribution $P(d\theta|x_n)$. 
We assume Bernestein von-Mises's theorem, that is, for some $c_n\rightarrow\infty$ and 
some $\mathbb{R}^d$-valued random variable $u_n(x_n)$ such that 
\begin{equation}\label{co}
c_n(\theta-u_n)=O_\mathbb{P}(1)
\end{equation}
where $\theta|x_n\sim P(d\theta|x_n)$. We will consider asymptotic properties of 
the scaled parameter $c_n(\theta-u_n)$. 

MCMC procedure generates a sequence $\theta_\infty=(\theta_0,\theta_1,\ldots)$ such that 
the law of $\theta_\infty|x_n$ is a Markov chain
with the invariant distribution $P(d\theta|x_n)$. 
We assume stationarity of the process $\theta_\infty|x_n$, that is, the initial guess $\theta_0$ is generated from 
the posterior distribution. This is impractical setting, but this assumption can be weakened. See \citet{Kamatani10} for the detail. 
Let 
\begin{equation}
I_M^n(f)=\frac{1}{M}\sum_{m=0}^{M-1}f(c_n(\theta_i-u_n)),\ I^n(f)=\int f(c_n(\theta-u_n))P(d\theta|x_n). 
\end{equation}
We expect that $I_M^n(f)$ is a good approximation of $I^n(f)$. For this in mind, we define local consistency. 

\begin{definition}[Local consistency]
MCMC procedure 
is said to have the local consistency if 
$I_{M_n}^n(f)-I^n(f)=o_{\mathbb{P}}(1)$ for any continuous,  bounded function $f$
and for any $M_n\rightarrow\infty$. 
\end{definition}

The MCMC procedure does not  always work well. 
We also define a property of this inefficient behavior. Essentially, the good behavior, local consistency, and the bad behavior, local degeneracy defined below are exclusive (see \citet{Kamatani11b}). 

\begin{definition}[Local degeneracy]
MCMC procedure 
is said to have the local degeneracy if 
$I_{M}^n(f)-I_1^n(f)=o_{\mathbb{P}}(1)$ for any continuous,  bounded function $f$
and for any $M\in\mathbb{N}$. 
\end{definition}

As a measure of poor behavior, local degeneracy is sometimes too wide
and so we define a kind of order of convergence among degenerate MCMC procedures. 

\begin{definition}[Local weak consistency]\label{wc}
MCMC procedure
is said to have the local weak consistency, if
$I_{M_n}^n(f)-I^n(f)=o_{\mathbb{P}}(1)$ for any continuous,  bounded function $f$
and for any $M_n$
such that $M_n/r_n\rightarrow \infty$. 
We call $r_n$ the order of the local weak consistency. 
\end{definition}

We can interpret $M_n$ as the sufficient number of iterations for good approximation. Therefore if $r_n$ is large, 
the MCMC procedure  requires many iterations to have a good result. 
Under the local consistency, we can take $r_n=1$. We can compare different algorithms by this order $r_n$.

\subsection{Useful lemma}\label{WC}
Let $r_n\rightarrow\infty$. 
For each $n=1,2,\ldots$, consider a $\mathbb{R}^d$-valued semi-Markov process $\{\theta_t^n;t\ge 0\}$, which 
jumps at $r_n^{-1}, 2r_n^{-1}, \ldots$
on a probability space $(\Omega_n,\mathcal{F}_n,\mathbb{P}_n)$. 

\begin{lemma}\label{PRAthe1}
Assume that the embedded Markov chain $\{\theta_{m/r_n}^n;m=0,1,\ldots\}$ is a stationary Markov chain with the invariant distribution $P_n$. 
If $\{\theta_t^n;t\ge 0\}$ converge in law to a stationary  ergodic process, then 
\begin{equation}\label{PRAeq0}
\frac{1}{T_n}\int_{[0,T_n)}f(\theta_t^n)dt-\int_{\mathbb{R}^d}f(\theta)P_n(d\theta)=o_{\mathbb{P}_n}(1)
\end{equation}
for any $T_n\rightarrow\infty$ for any bounded and continuous function $f$. 
\end{lemma}

\begin{proof}
%Without loss of generality, we can assume 
%Without loss of generality, we may assume 
% $P_n(M^n_x\notin\mathcal{S}_{ c_n})=0$.
We omit the subscript $n$ from $\theta_t^n$. 
Proof is almost the same as that of Lemma 2 of \citep{Kamatani10}. 
Without loss of generality, we can assume $\sup_{x\in \mathbb{R}^d}|f(x)|\le 1$. 
Write $I_{T_n}$ and $I^n$ for the first and the second terms in the left-hand side of (\ref{PRAeq0}), respectively,
and write $I_{i,S}$ for 
$\frac{1}{S}\int_{t\in [0,S)}f(\theta_{iS+t})dt$.  Then 
\begin{equation}\nonumber
I_T=\frac{S}{T}\sum_{i=0}^{[T/S]-1}I_{i,S}+\frac{1}{T}\int_{[S[T/S],T)}f(\theta_t)dt. 
\end{equation}
Note that $I_{i,S}\ (i=0,1,\ldots)$ is not identically distributed in general. If we take
 $S_n=([Sr_n]+1)/r_n$, then $I_{i,S_n}\ (i=0,1,\ldots)$ have the same law under $\mathbb{P}_n$. 
Hence as in Lemma 2 of \citep{Kamatani10}, we have
\begin{eqnarray*}\nonumber
\mathbb{E}_n[|I_T-I^n|]&\le& \frac{S_n}{T}\Big[\frac{T}{S_n}\Big]\mathbb{E}_n[|I_{0,S_n}-I^n|]+2\frac{T-S_n[T/S_n]}{T}\\
&\le & \frac{S_n}{T}\Big[\frac{T}{S_n}\Big]\Big\{\frac{S}{S_n}\mathbb{E}_n[|I_{0,S}-I^n|]+2\frac{S_n-S}{S_n}\Big\}+2\frac{T-S_n[T/S_n]}{T}
\end{eqnarray*}
where in the second inequality, we used $I_{0,S_n}=(S/S_n)I_{0,S}+\int_{S}^{S_n}f(\theta_t)dt/S_n$. 
As $T=T_n\rightarrow\infty$, the second and the  third terms in the left-hand side vanishes, and $S_n/T[T/S_n]S/S_n\rightarrow 1$ for any fixed $S$. 
Hence if $\mathbb{E}_n[|I_{0,S}-I^n|]$ can be arbitrary small, the claim follows. 

Let $P$ be the limit of $P_n$
and write $I=\int_{\mathbb{R}^d}f(\theta)P(d\theta)$. Then 
$\mathbb{E}_n[|I_{0,S}-I^n|]\rightarrow \mathbb{E}[|I_{0,S}-I|]$
where $\mathbb{E}$ is the expectation with respect to the limit probability measure $\mathbb{P}$.
Since $\theta_t$ is ergodic under $\mathbb{P}$, this value tends to $0$ as $S\rightarrow\infty$. Hence the claim follows. 
\end{proof}

\section{Application to mixture model}\label{GSM}
Let $F_t(dx)=f_t(x)dx$ be probability measures on $(E,\mathcal{E})$
with parameter $t\ge 0$, and write $F=F_0$ and $f=f_0$. 
We assume that $f_t$ is always strictly positive.  
Consider the following simple mixture model:
\begin{equation}\label{simplemixture}
P^\epsilon(dx|\theta)=p^\epsilon_\theta(x)dx=(1-\theta)F(dx)+\theta F_\epsilon(dx).
\end{equation} 
MCMC procedures for general $k$-component mixture model have been developed to perform better posterior inference. 
See monographs such as \citet{RC} and  \citet{opac-b1129057}. 
It is well known that for general $k$-component mixture model, the posterior distribution is multi-modal, and if these peeks are close, then
the posterior inference becomes difficult due to the so-called label-switching problem (see \citet{SJRSSB00,Marin2005459} and \citet{Jasra2005}). 
We address here a separate issue. 
In fact, under such a situation, another problem, local degeneracy occurs. 
We illustrate this effect by using the order of the weak consistency. 

For this reason, we  assume over-parametrized situation, that is, 
the observation $x_n=\{x^1,\ldots, x^n\}$ are independent draw from a one-component model $F$. 
We will show that if two components $F_0$ and $F_\epsilon$ are close, the performance becomes even worse. 
To illustrate the effect, 
we let $\epsilon=\epsilon_n\rightarrow 0$. 
Write $r_n=\epsilon^{-1}n^{1/2}$ and $ c_n=\epsilon n^{1/2}$.
There is an obvious relation $r_nc_n=n$. 
As in p902 of \citet{Gassiat2002897}
we assume the following regularity condition. 
Write $L^2(F)$ for the set of $F$-square integrable functions with norm $\|\cdot\|$ defined by 
$\|g\|^2=\int g(x)^2 F(dx)$. 

\begin{assumption}\label{Assumption}
There exists $d\in L^2(F)$ such that
$r_\epsilon(x):=f_\epsilon(x)/f(x)-1-\epsilon d(x)=o(\epsilon)$ in $L^2(F)$. 
Moreover, 
$\|d\|^2=I\neq 0$. 
The prior distribution is assumed to be $\mathrm{Beta}(\alpha_1,\alpha_0)$ for $\alpha_0, \alpha_1>0$. 
\end{assumption}

Under the assumption, $\int d(x)F(dx)=0$ because $|\int d(x)F(dx)|=\epsilon^{-1}|\int r_\epsilon(x)F(dx)|\le \epsilon^{-1}\|r_\epsilon\|\rightarrow 0$. 
One step of the DA procedure is
\begin{equation}\label{damix}
\left\{\begin{array}{l}
y^i|x^i,\theta\sim \mathrm{Bernoulli}(p_i),\ p_i=\frac{\theta f_\epsilon(x^i)}{(1-\theta)f(x^i)+\theta f_\epsilon(x^i)}\ (i=1,\ldots, n),\\
\theta|x_n,y_n\sim \mathrm{Beta}(\alpha_1+n_1,\alpha_0+n-n_1)
\end{array}\right.
\end{equation}
where $y_n=\{y^1,\ldots, y^n\}$ and $n_1$ is the number of heads in $y_1,\ldots, y_n$. 
For this model, it is natural to take state space scaling as
\begin{equation}\label{ss}
\theta\mapsto c_n\theta
\end{equation}
where we take $u_n\equiv 0$ in (\ref{co}). We will discuss the local consistency and the local degeneracy 
under this localization. 

The DA output $\theta_0,\theta_1,\ldots$ behaves poorly, and the sequence converges to the following diffusion process after the state space scaling in (\ref{ss}) with suitable time scaling. 
Let
\begin{equation}\label{CIRlike}
dX_t=(\alpha_1+X_tz-X_t^2I)dt +\sqrt{2X_t}dW_t;\ X_0|z\sim P^*(dx|z)\propto \exp(x z-x^2I/2)x^{\alpha_1-1}dx
\end{equation}
where $\{W_t;t\ge 0\}$ is the standard Brownian motion independent of $z\sim N(0,I)$. 
Write the law of $\{X_t;t\ge 0\}$ given $z$ by $\mathcal{L}(\{X_t;t\ge 0\}|z)$. 
For each $z$, there is a weak solution $\mathcal{L}(\{X_t;t\ge 0\}|z)$ which is  ergodic with invariant measure $P^*(dx|z)$ (see Theorem 2.3 in \citet{MR2132002}. See also Section 5.5 of \citet{Karatzas1991}). 
By  convergence to this process, we obtain the following. 

\begin{theorem}\label{gibbsweak}
Suppose that $x^1,\ldots, x^n$ are independent sample from $F$. 
Under Assumption \ref{Assumption}, 
the DA procedure has the local weak consistency with the order $r_n=\epsilon^{-1}n^{1/2}$ if $c_n=\epsilon n^{1/2}\rightarrow\infty$.
\end{theorem}

\begin{proof}
Let $\{\theta_m^n;m=0,1,\ldots\}$ be the stationary Markov chain generated by (\ref{damix}). 
Let $X_t^n=c_n\theta_{[r_n t]}^n$.
Then by Theorem \ref{maincor}, 
$\mathcal{L}(\{X_t^n;t\ge 0\}|x_n)$ 
tends  to $\mathcal{L}(\{X_t;t\ge 0\}|z)$ in distribution, where $z\sim N(0,I)$. 
As mentioned above, $\mathcal{L}(\{X_t;t\ge 0\}|z)$ is stationary and ergodic. 
Together with the separability of the Skorohod topology, and Skorohod's representation theorem (see Theorem 6.7 of \citet{MR1700749}), 
there is a probability space such that 
$\mathcal{L}(\{X_t^n;t\ge 0\}|x_n)\rightarrow \mathcal{L}(\{X_t;t\ge 0\}|z)$. 
Hence by Lemma \ref{PRAthe1}, for any bounded continuous function $f$ on $[0,\infty)$, 
\begin{equation}\nonumber
\frac{1}{T_n}\int_{[0,T_n)}f(X_t^n)dt-
\int_{[0,1]} f(c_n\theta)P(d\theta|x_n)
=o_{\mathbb{P}_n}(1)
\end{equation}
for any $T_n\rightarrow\infty$, where $P(d\theta|x_n)$ is the posterior distribution, which is the invariant distribution 
of $X_t^n$. 
Take $M_n/r_n=T_n$ and rewrite 
\begin{equation}\nonumber
\frac{1}{T_n}\int_{[0,T_n)}f(X_t^n)dt=\frac{1}{M_n}\sum_{m=0}^{M_n-1}f(c_n\theta_m).
\end{equation}
Then the convergence of probability in the above means weak consistency of the DA procedure on the order $r_n$. 
\end{proof}

Although this result for the large sample scaling limit is new, the scaling limit to a diffusion process have been studied in other directions by \citet{gmjota91} for small variance asymptotics, and  by \citet{MR1428751}
for high-dimensional small variance asymptotics. In particular, the latter approach is still very active. See a recent paper by \citet{MR2977981} and the references their in. 
It is worth mentioning that  we can apply the local weak consistency for these results.

\section{Numerical results}\label{NR}

\subsection{Metropolis-Hastings procedure}
To illustrate poor performance of the DA procedure, we consider a simple independent type Metropolis-Hastings (IMH) procedure
as an alternative and compare it with the DA procedure. Note that we prepare this IMH procedure just for comparison  and may work well only for this simple mixture model. However related methods may work well for general mixture model and this direction will be studied in elsewhere. 

We briefly review the IMH procedure. 
If we want to approximate probability distribution $F$ on $(E,\mathcal{E})$, we prepare the so-called proposal distribution $G$, such that there exists a Radon-Nikodym derivative $dF/dG=h(x)$. Then IMH procedure iterates the following;
Suppose that we have the current value $\theta\in E$. Then 
\begin{equation}\nonumber
\mathbf{simulate}\ \theta^*\sim G,\ \mathbf{set}\ \theta\leftarrow\left\{\begin{array}{ll}\theta^*\ \mathrm{with\ probability}\ \alpha(\theta,\theta^*)\\
\theta\ \mathrm{with\ probability}\ 1-\alpha(\theta,\theta^*)\end{array}\right.,
\end{equation}
where $\alpha(\theta,\theta^*)=\min\{1, h(\theta^*)/h(\theta)\}$. 
This iteration resulted in a Markov chain $\theta_0,\theta_1,\ldots$ with the invariant distribution $F$. Hence if it is ergodic, we obtain an approximation of $F$ without simulation from $F$. 

Now we apply this IMH procedure to the simple mixture model. 
The key is the choice of the proposal distribution. 
Set $Q^\epsilon(dx|\theta)$ which is close to $P^\epsilon(dx|\theta)$ in such as the Kullback-Leibler distance. Calculate the posterior distribution $Q(d\theta|x_n)$
for observation $x_n=\{x^1,\ldots, x^n\}$ with the model $Q(dx|\theta)$. 
We use $Q^\epsilon(d\theta|x_n)$ as the proposal distribution. 

Next section, we will consider $F_t=N(t,1)$. For this, take $Q^\epsilon(dx|\theta)=N(\epsilon\theta,1)$
with the uniform prior. Then 
$Q(d\theta|x_n)=N(\sum_{i=1}^n x^i/n\epsilon,1/n\epsilon^2)$ truncated to $[0,1]$. 
It is not difficult to check the local consistency of this IMH procedure, but it is beyond our scope. 

\subsection{Simulation}
\begin{figure}
   \includegraphics[width=1\textwidth,bb=0 0 864 432]{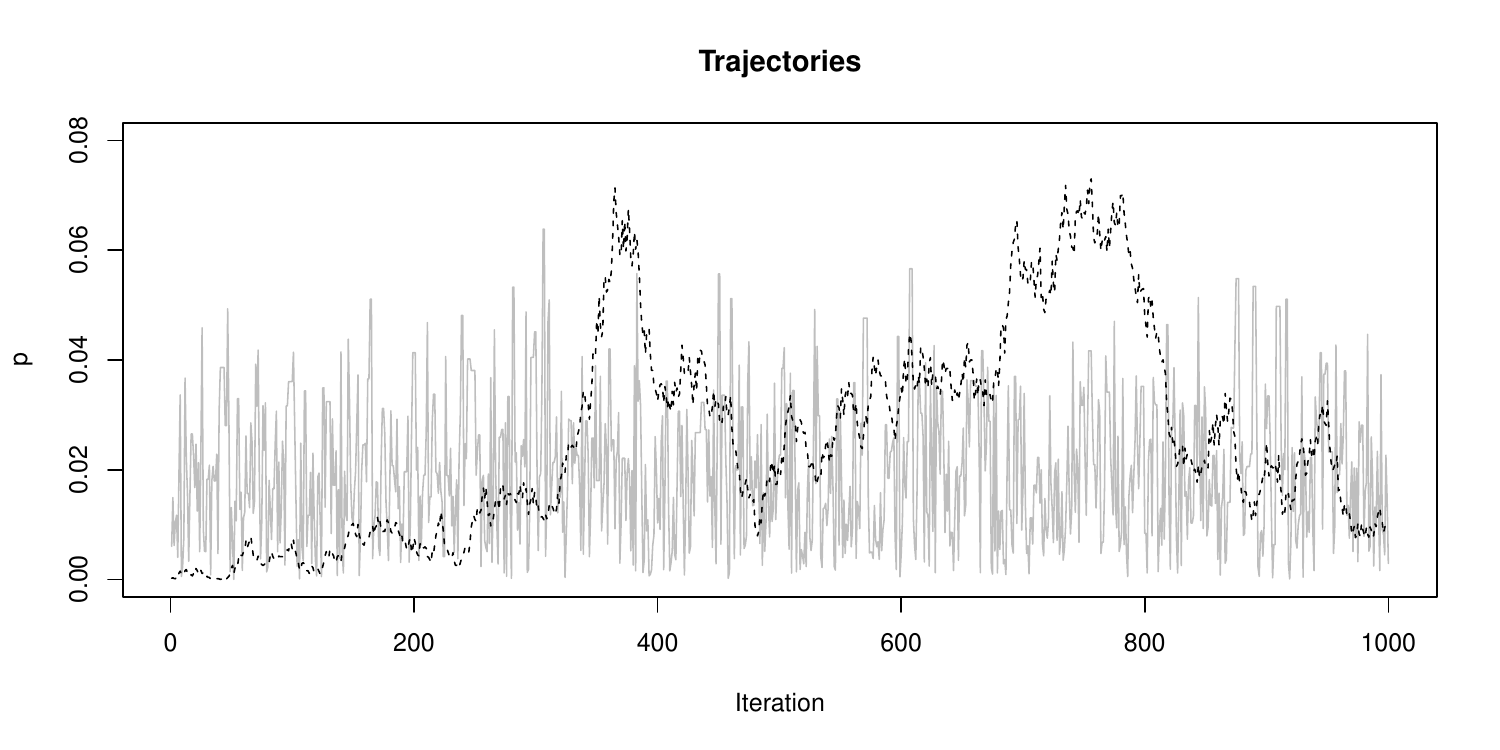}
   \caption{Trajectories of MCMC procedure for $n=10^4$. The dashed line is the trajectory from the DA procedure and 
the solid gray line is that of the IMH procedure. }
  \label{fig:two}
\end{figure}

We compare the DA and the IMH  procedures through  numerical 
simulations. 
Consider the normal mixture model $F_\epsilon(dx)=N(\epsilon,1)$. 
To illustrate the difference of the DA and the IMH procedures, first we plot the trajectories of $\theta_m\ (m=0,1,2,\ldots)$
under fairly large sample size $n=10^4$ with the true model $N(0,1)$ and $\epsilon=0.5$. 
Unlike the IMH procedure, the trajectory from the DA procedure behaves like a stochastic diffusion process (Figure \ref{fig:two})
and this is true by Theorem \ref{maincor}. 
By Theorem \ref{gibbsweak}, the order of the weak consistency is $r_n=\epsilon^{-1}n^{1/2}=200$ for the DA procedure  but it is $r_n=1$ for the IMH procedure. 
 
\begin{figure}[htbp]
\includegraphics[width=5cm,bb=0 0 504 504]{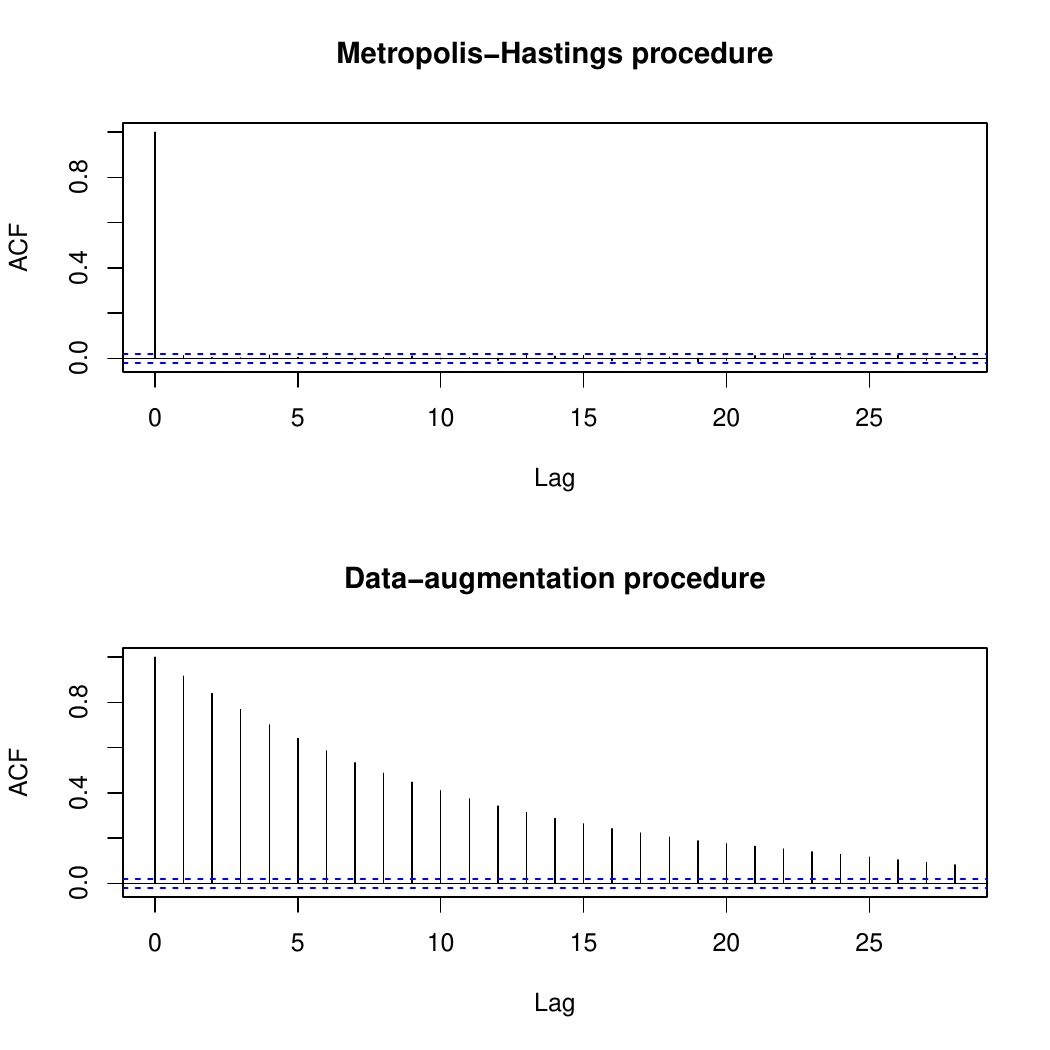}
\includegraphics[width=5cm,bb=0 0 504 504]{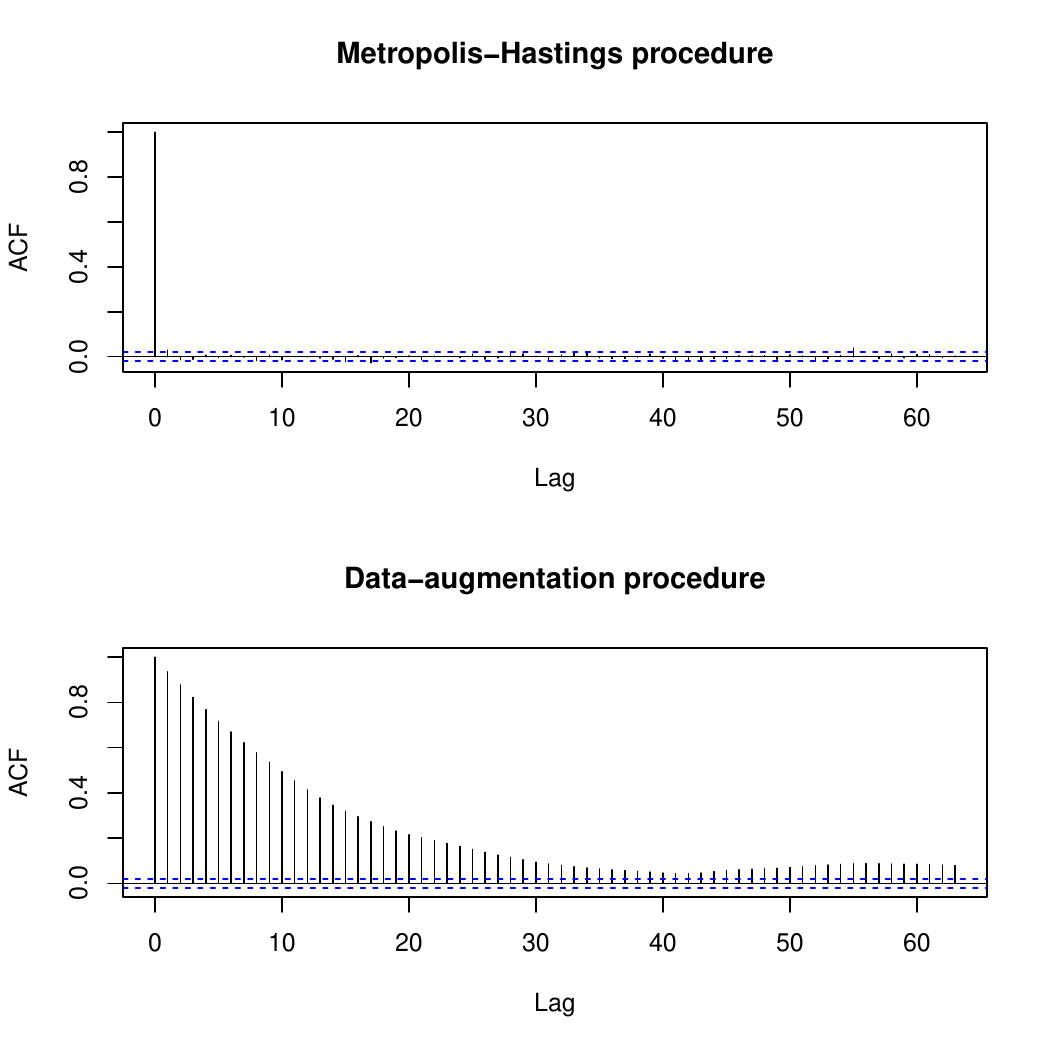}
\includegraphics[width=5cm,bb=0 0 504 504]{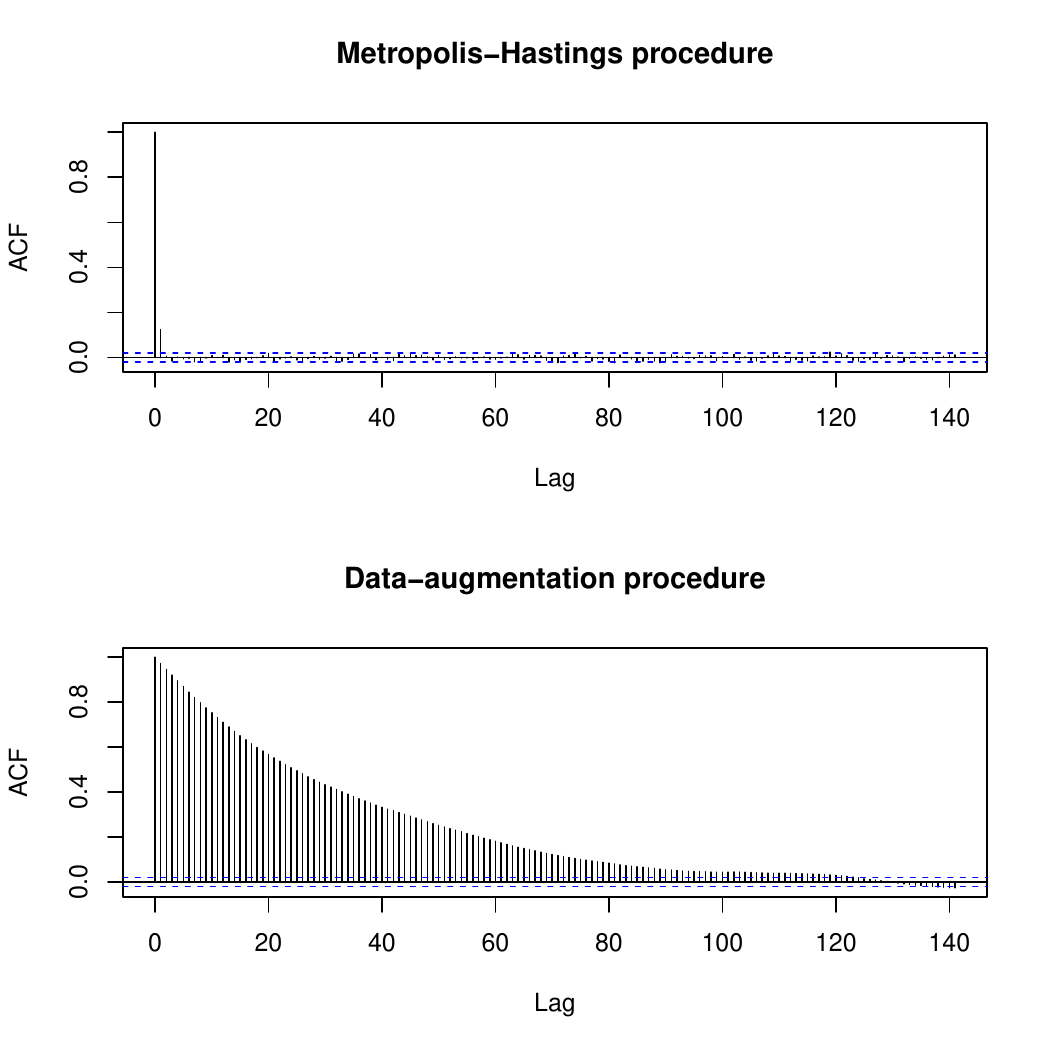}
\caption{Empirical autocorrelations for $n=50$ (left), $250$ (middle), and $1250$ (right), 
for the IMH procedure (top) and the DA procedure (bottom). }
  \label{fig:three}
\end{figure}

Next we check the effect of $n$, $\epsilon$ and the underlying true model. To illustrate the differences of the performance, 
we plot empirical autocorrelations. If an MCMC procedure has poor mixing property, empirical autocorrelation
does not converges to $0$ quickly. 
First we check the effect of $n$ for $\epsilon=0.5$ by the different sample sizes $50, 250$ and $1250$. 
Orders of the weak consistency of the DA procedure are $r_n=14.14\ldots, 31.62\ldots$ and $70.71\ldots$. Recall that $r_n$ corresponds to the 
number of iteration for good convergence, and so we take the window size as $\max\{25,2r_n\}$. 
As the sample size becomes larger,  the mixing property of the DA procedure becomes worse, as the empirical autocorrelations suggest (Figure \ref{fig:three}). 

\begin{figure}[htbp]
\includegraphics[width=5cm,bb=0 0 504 504]{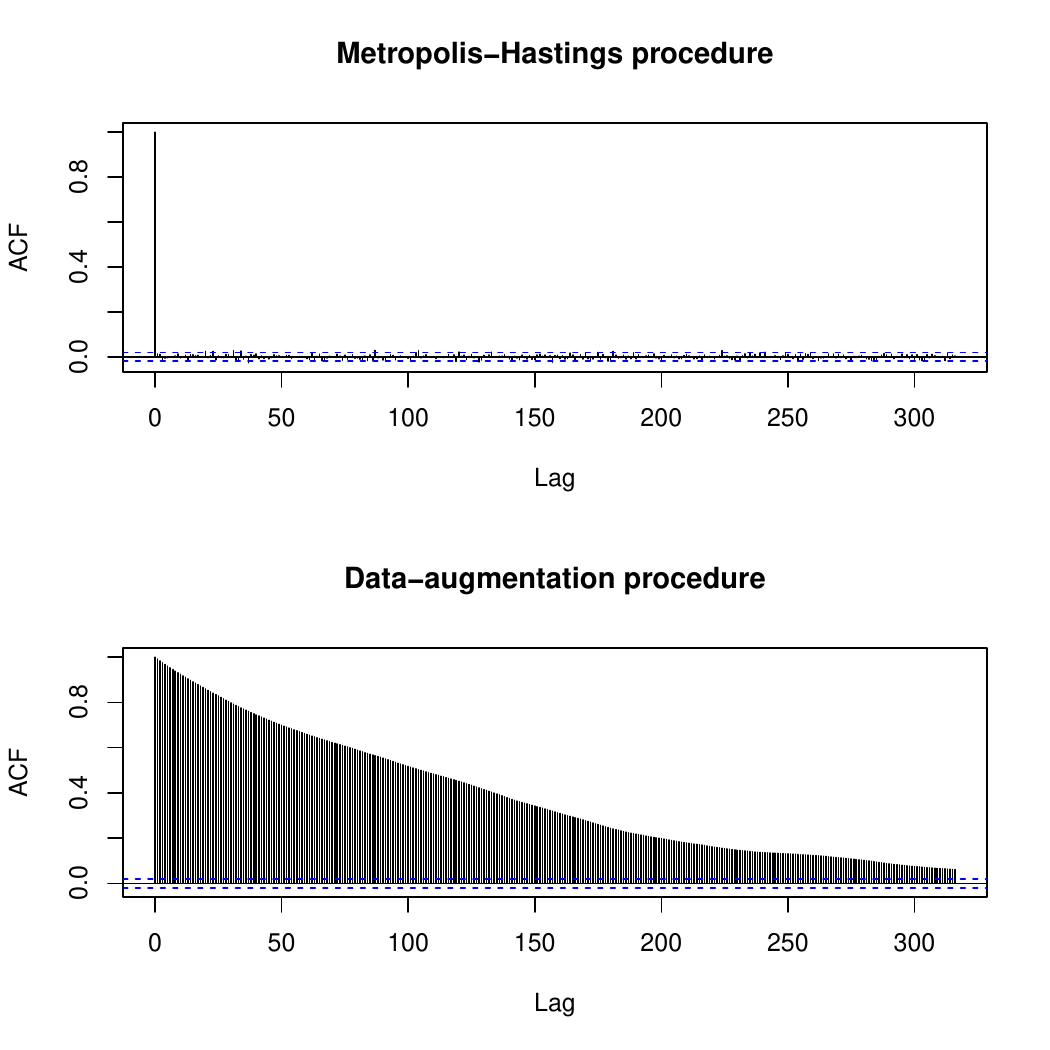}
\includegraphics[width=5cm,bb=0 0 504 504]{mixture_fixt_mh_acf-250-5-0.pdf}
\includegraphics[width=5cm,bb=0 0 504 504]{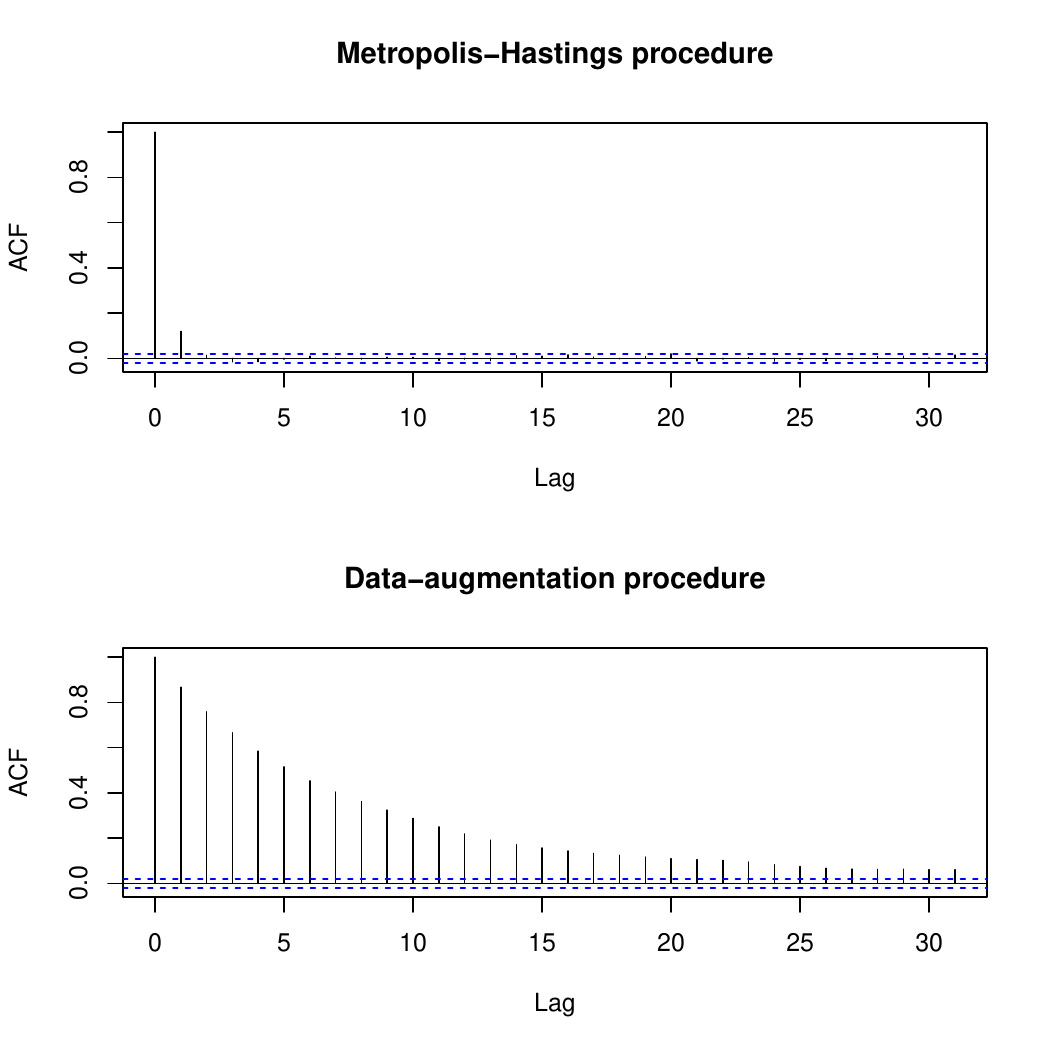}
\caption{Empirical autocorrelations for $\epsilon=0.1$ (left), $0.5$ (middle), and $1$ (right) with the sample size $n=250$, 
for the IMH procedure (top) and the DA procedure (bottom). }
  \label{fig:four}
\end{figure}

Similarly, next we check the effect of $\epsilon$ for $n=250$ by the different values $\epsilon=0.1, 0.5$ and $1$. 
Orders of the weak consistency of the DA procedure are $r_n=158.11\ldots, 31.62\ldots$ and $15.81\ldots$ (Figure \ref{fig:four}). 

\begin{figure}[htbp]
\includegraphics[width=5cm,bb=0 0 504 504]{mixture_fixt_mh_acf-250-10-0.pdf}
\includegraphics[width=5cm,bb=0 0 504 504]{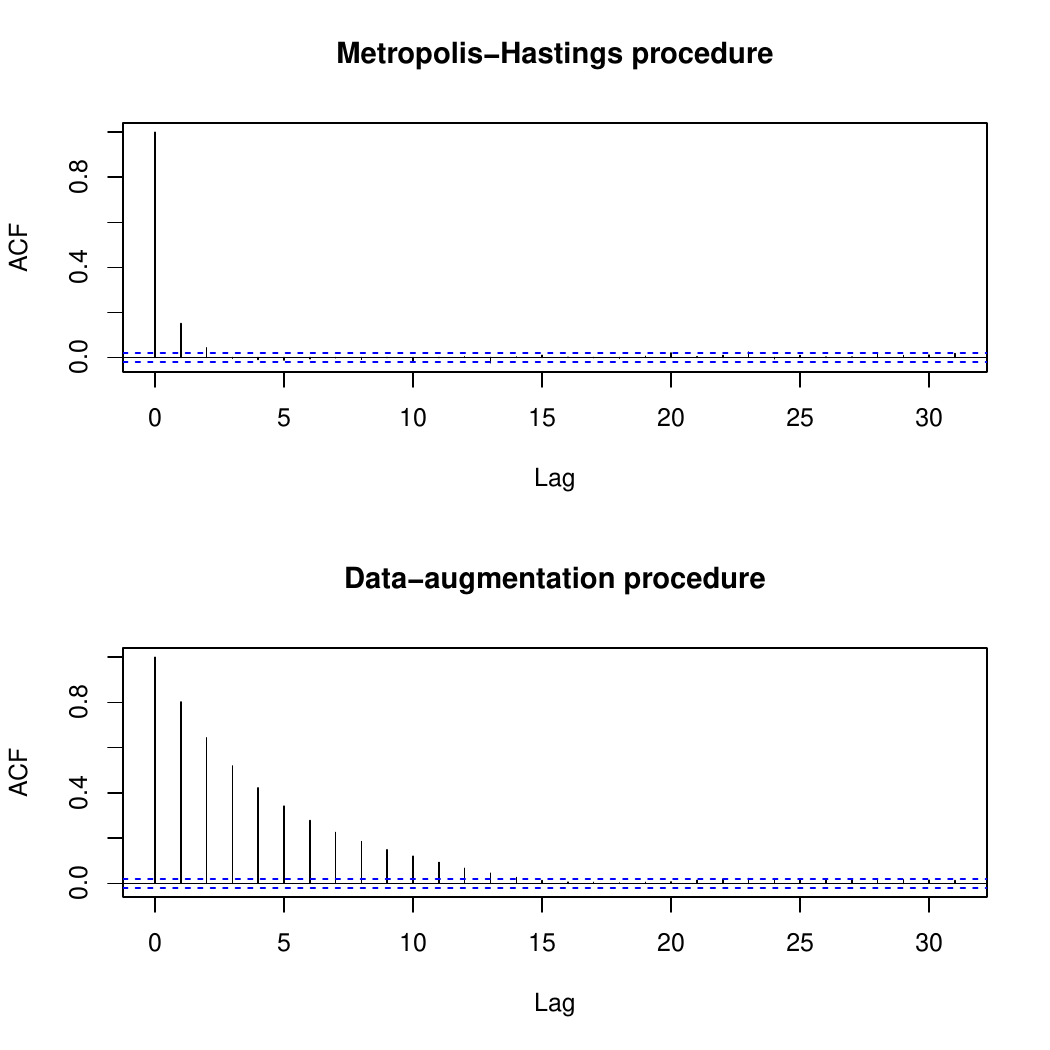}
\includegraphics[width=5cm,bb=0 0 504 504]{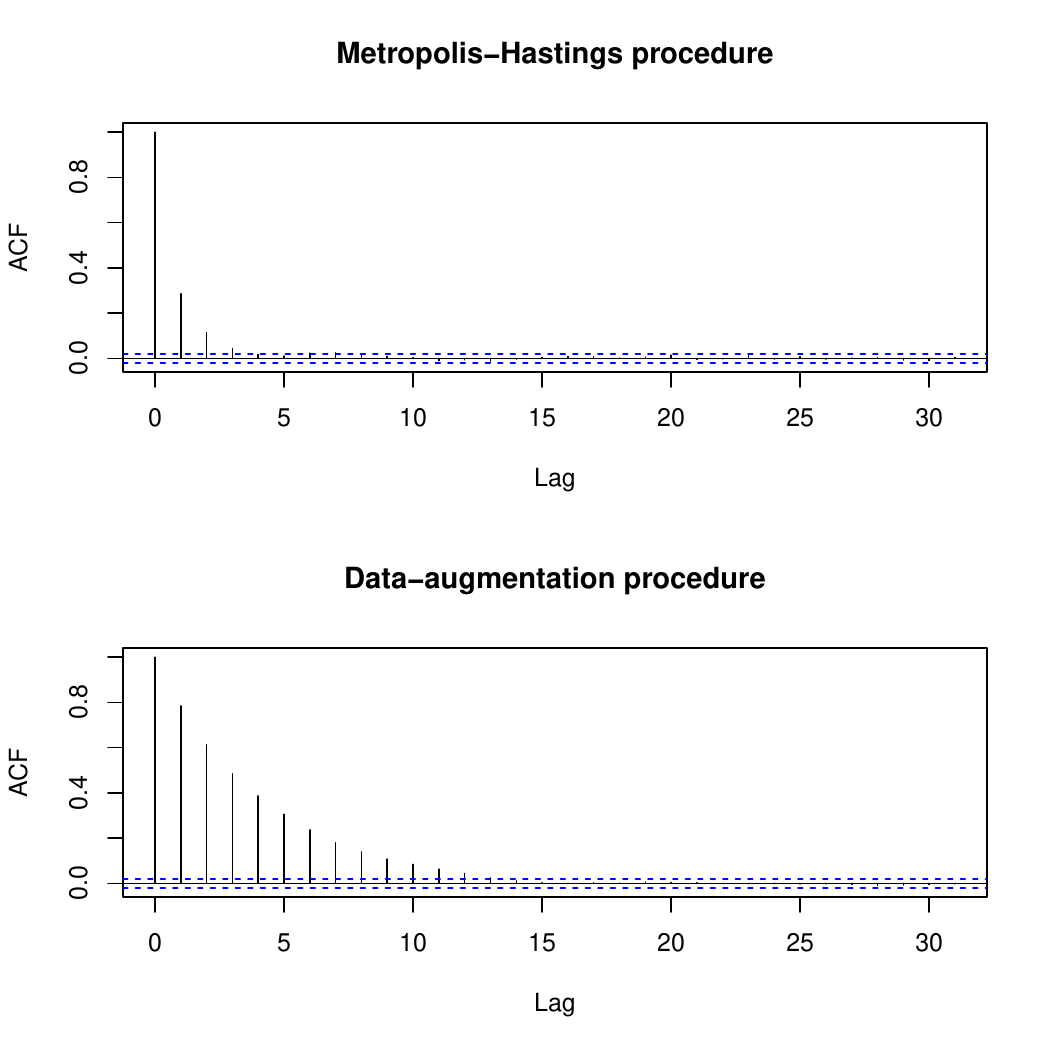}
\caption{Empirical autocorrelations for $p=0$ (left), $0.25$ (middle), and $0.5$ (right) with the sample size $n=250$, 
and $\epsilon=1$
for the IMH procedure (top) and the DA procedure (bottom). }
  \label{fig:five}
\end{figure}

Finally, we check the effect of the underlying model. We only checked the behaviors of MCMC procedures under the true model 
$N(0,1)$. Now we check those for $pN(0,1)+(1-p)N(1,1)$ for $p=0$, $0.25$ and $0.5$. 
However as Lemma \ref{lem1} and Le Cam's third lemma suggest (see the comment after Lemma \ref{lem1}), 
the effect of the difference of the underlying true model
may be small under the assumption $\epsilon n^{1/2}\rightarrow\infty$ (Figure \ref{fig:five}).  

%\footnote{specify epsilon for Figure 2}. 

\subsection{Discussion}
In this article we presented a definition of the order of the weak consistency, 
and applied it to a simple mixture model. The simulation results suggest that this order is a good measure of the 
convergence. It remains true that the verification of the weak consistency for complicated model is technically difficult. 
However estimation of $r_n$ is not difficult if the corresponding Markov chain convergence to a diffusion process such as (\ref{sde}). 
Such a convergence is probably true for more general MCMC methods. For  example, the DA procedure for a simple probit model converges to the Ornstein-Uhlenbeck process with the rate $r_n=n$.  With this in mind, $r_n$ can be estimated by
an empirical estimate of 
\begin{equation}\nonumber
\frac{\mathbb{V}_n[\theta|x_n]}{\mathbb{V}_n[\theta_1^n-\theta_0^n|x_n]}\approx r_n
\end{equation}
where the numerator is the variance of the posterior distribution, and the denominator is the variance of 
$\theta_1^n-\theta_0^n$, where $\{\theta_m^n;m=0,1,\ldots\}$ is an output of an MCMC procedure. 
A similar statistic was proposed in \citet{RSSB:RSSB123}. 
This can be used as a measure of efficiency of MCMC methods.

\section*{Acknowledgements}
I am deeply grateful to Professor Nakahiro Yoshida, for encouragement and valuable suggestions.

\appendix 
\section{Appendix}

\subsection{Some properties of simple mixture model}\label{INTmix}
We define 
\begin{equation}\nonumber
t_g=\frac{g-f}{f},\ 
l_n(g)=\sum_{i=1}^n\log g(x^i). 
\end{equation}　
For $g=p_\theta^\epsilon$, we have $t_{p^\epsilon_\theta}=\theta t_{f_\epsilon}$
and $t_{f_\epsilon}=\epsilon d+r_\epsilon$. 
Note that $\int g(x)F(dx)=\int r_\epsilon(x)F(dx)=0$. 

\begin{lemma}[Local asymptotic normality]\label{lem1}
Under Assumption \ref{Assumption},
for any $H>0$, 
\begin{equation}\label{unifconv}
\sup_{h\in [0,H]}|l_n(p^\epsilon_{ h/c_n})-l_n(f)-hn^{-1/2}\sum_{i=1}^n d(x^i)+h^2I/2|=o_{\mathbb{P}_n}(1).
\end{equation}
\end{lemma}

\begin{proof}
Set $R$ to be  $\log(1+x)=x-x^2/2+x^2R(x)$, where $R(x)\rightarrow 0\ (x\rightarrow 0)$. 
As in  the proof of Theorem 3.1 of \citet{Gassiat2002897} we have
\begin{displaymath}
l_n(g)-l_n(f)=\sum_{i=1}^n t_g(x^i)-\frac{1}{2}\sum_{i=1}^nt_g(x^i)^2+\sum_{i=1}^n t_g(x^i)^2R(t_g(x^i)).
\end{displaymath}
For $g=p_{h/c_n}^\epsilon$, we have $t_g=hn^{-1/2}(d+r_\epsilon/\epsilon)$. Hence the first term of the right-hand side is 
\begin{equation}\nonumber
h n^{-1/2}\sum_{i=1}^n d(x^i)+hn^{-1/2}\epsilon^{-1}\sum_{i=1}^n r_\epsilon(x^i)=h n^{-1/2}\sum_{i=1}^n d(x^i)+o_{\mathbb{P}_n}(1)
\end{equation}
uniformly in $h\in [0,H]$. 
By similar arguments, $\sup_{h\in [0,H]}|\sum_{i=1}^nt_{p_{h/c_n}^\epsilon}(x^i)^2-h^2 I|=o_{\mathbb{P}_n}(1)$ in probability. For the third term, 
\begin{equation}\nonumber
\sup_{h\in [0,H]}\max_{i=1,\ldots, n}| t_{p^\epsilon_{h/c_n}}(x^i)|
\le 
H\max_{i=1,\ldots, n}|c_n^{-1} t_{f_\epsilon}(x^i)|. 
\end{equation}
However 
$\max_{i=1,\ldots, n}|c_n^{-1} t_{f_\epsilon}(x^i)|=o_{\mathbb{P}_n}(1)$
by the inequality $\mathbb{P}(\max_{i=1,\ldots, n}|X_i|>c)\le c^{-2}\sum_{i=1}^n\mathbb{E}[|X_i|^2,|X_i|>c)]$. Hence
the third term is negligible. 
 By these, the convergence (\ref{unifconv}) follows. 
%$P_n$ and $P_n_{s_n}$ are mutually contiguous. 
%
%Latter claim is easy to verify since $P_n(\psi_n)\rightarrow 0$ for $\psi_n:X_n\rightarrow [0,1]$, then 
%$p(\psi_n|\theta= h/c_n)\rightarrow 0$ for any $h\in [0,H]$, which yields $Q_n(\psi_n)\rightarrow 0$. 
%On the other hand, if $Q_n(\psi_n)\rightarrow 0$, then 
%$p(\psi_n|\theta= h/c_n)\rightarrow 0$ for a.s. $h$, and hence we can choose convergent subsequence 
%$s_n\rightarrow h>0$ among then. Then by the former claim, $P_n(\psi_n)\rightarrow 0$. 
\end{proof}

Let $P_n(dx_n)=\prod_{i=1}^n F_0(x^i)$
and
$Q_n(dx_n)=\prod_{i=1}^n P^\epsilon_{h/c_n}(x^i)$ where $x_n=\{x^1,\ldots, x^n\}$. Then $P_n$ and $Q_n$ are mutually contiguous by Le Cam's third lemma. 
Therefore if a statement is true in $P_n$-probability, then it is also true in $Q_n$-probability, and vice versa.
Thus Theorem \ref{gibbsweak} holds even if the observation is an i.i.d. draw from $P^\epsilon_{h/c_n}$ for any fixed $h>0$.

\begin{lemma}[Consistency of the posterior distribution]\label{lardevi}
Under Assumption \ref{Assumption},
for any $M_n\rightarrow \infty$, 
\begin{displaymath}
\int_{M_n/c_n}^1P(d\theta|x_n)=o_{\mathbb{P}_n}(1). 
\end{displaymath}
\end{lemma}

\begin{proof}
Let $P(d\theta)=\mathrm{Beta}(\alpha_1,\alpha_0)$. By definition, for $M, H>0$, 
\begin{eqnarray}
\int_{M/c_n}^1P(d\theta|x_n)&=&\frac{\int_{M/c_n}^1\exp(l_n(p_\theta^\epsilon)-l_n(f))P(d\theta)}{\int_0^1\exp(l_n(p_\theta^\epsilon)-l_n(f))P(d\theta)}\le 
\frac{\int_{M/c_n}^1\exp(l_n(p_\theta^\epsilon)-l_n(f))P(d\theta)}{\int_0^{H/c_n}\exp(l_n(p_\theta^\epsilon)-l_n(f))P(d\theta)}\nonumber\\
&\le & 
\frac{\int_{M/c_n}^1\exp(l_n(p_\theta^\epsilon)-l_n(f))P(d\theta)}{\inf_{\theta\in [0,H/c_n]}\exp(l_n(p_\theta^\epsilon)-l_n(f))\int_0^{H/c_n}P(d\theta)}\label{appendix_estim}
\end{eqnarray}
 By Lemma \ref{lem1}, the infimum in the denominator is bounded below from $0$ in probability, 
and 
$\int_0^{H/c_n}P(d\theta)$ is on the order of $c_n^{-\alpha_1}$. 
Therefore, to prove the claim, it is sufficient to show that the numerator is $o_{\mathbb{P}_n}(c_n^{-\alpha_1})$
as $M=M_n\rightarrow\infty$. 

As in the proof of Inequality 1.2 of \citet{Gassiat2002897}, since $\log(1+x)\le x-(x_-)^2/2$
where $x_-=\min\{0,x\}$, we have
\begin{displaymath}
l_n(p_\theta^\epsilon)-l_n(f)\le \theta\sum_{i=1}^n t_{f_\epsilon}(x^i)-\frac{\theta^2}{2}\sum_{i=1}^n (t_{f_\epsilon})^2_-(x^i).
\end{displaymath}
Note that $|(t_{f_\epsilon})_--\epsilon d_-|\le |r_\epsilon|$. 
Since $P(d\theta)\le C\theta^{\alpha_1-1}d\theta$ for some constant $C>0$, 
up to this constant,  an upper bound of the numerator of (\ref{appendix_estim}) is
\begin{eqnarray*}\nonumber
\int_{M/c_n}^\infty \exp(\theta\sum_{i=1}^n t_{f_\epsilon}(x^i)-\frac{\theta^2}{2}\sum_{i=1}^n (t_{f_\epsilon})^2_-(x^i))\theta^{\alpha_1-1}d\theta&=:&c_n^{-\alpha_1}
\int_{M}^\infty \exp(hA_n-\frac{h^2}{2}B_n)
h^{\alpha_1-1}dh
\end{eqnarray*}
where $A_n:=c_n^{-1}t_{f_\epsilon}(x^i)\Rightarrow N(0,I)$
and $B_n:=\sum_{i=1}^n c_n^{-2}(t_{f_\epsilon})^2_-(x^i)\rightarrow \|d_-\|^2$ in probability. 
Hence as $M=M_n\rightarrow\infty$, this is on the order of $o_{\mathbb{P}_n}(c_n^{-\alpha_1})$ and hence the claim follows. 
\end{proof}

Let $P^*(d\theta|z)$ be a probability measure on $[0,\infty)$ such that
\begin{equation}\label{asymppost}
P^*(d\theta|z)\propto \exp(\theta z-\theta^2I/2)\theta^{\alpha_1-1}d\theta
\end{equation}
where $z\in\mathbb{R}$. 
The following is a Bernstein-von Mises theorem. We omit this proof since it come from exact the same way as the proof of the usual Bernstein-von Mises theorem for regular parametric family. 
See for example, p142 of \citet{V}. In the following, $\|\cdot\|$ is the total variation distance defined in (\ref{tv}). 

\begin{proposition}[Bernstein-von Mises theorem]\label{BvM}
Under Assumption \ref{Assumption}, 
\begin{displaymath}
\|P^*(d\theta|x_n)-P^*(d\theta|z=n^{-1/2}\sum_{i=1}^n d(x^i))\|=o_{\mathbb{P}_n}(1)
\end{displaymath}
where 
$P^*(d\theta|x_n)$ is defined by $\int_{[0,H]}P^*(d\theta|x_n)=\int_{[0,H/c_n]} P(d\theta|x_n)$. 
\end{proposition}

%\begin{proof}
%Let
%\begin{displaymath}
%r_{n,h}(x_n)
%=
%\log(L_{n,h})-hZ_n+h^2I/2+(\alpha_0-1)\log(1- h/c_n).
%\end{displaymath}
%By  Lemma \ref{lem1}, for any $\delta>0$, 
%setting $A_n=\{x_n;\sup_{h\in [0,H]}|r_n^h(x_n)|<\delta\}$, 
%we obtain $P_n(A_n)\rightarrow 1$. 
%By representation (\ref{eq4}) with the fact
%\[
%L_{n,h}p_\Theta^*(dh)=C\exp(hZ_n-h^2I/2+r_{n,h})h^{\alpha_1-1}dh
%\]
%for some constant $C>0$,  we obtain 
%\begin{displaymath}
%1_{[0,H]}(h)P^*(dh|x_n)\le P^*(dh|d=Z_n)\exp(2\delta)(1-P^*([0,H]^c|d=Z_n))^{-1}
%\end{displaymath}
%for $x\in A_n$. Therefore, 
%\begin{eqnarray*}
%\|P^*(\cdot|x_n)-P^*(\cdot|d=Z_n)\|&\le& 
%P^*([0,H]^c|x_n)+2(P^*(\cdot|x_n)-P^*(\cdot|d=Z_n))^+([0,H])\\
%&\le& P^*([0,H]^c|x_n)+2(\frac{\exp(2\delta)}{1-P^*([0,H]^c|d=Z_n))}-1)
%\end{eqnarray*}
%for $x\in A_n$. 
%Since 
%$P^*([0,H]^c|d=Z_n)\rightarrow 0$
%in $P_n$-probability, 
%taking $H\rightarrow\infty$ and $\delta\rightarrow 0$, 
%by Lemma \ref{lardevi}, the claim follows. 
%\end{proof}

\subsection{Convergence to a diffusion process}
Let $\theta^n_0, \theta^n_1,\ldots$ be an output of the DA procedure for given observation $x_n$. 
In this section, we will show that $\{h_m^n;m=0,\ldots,\}$ 
has a diffusion limit
where $h^n_m=c_n\theta^n_m$. 
%Write
%\begin{equation}\label{Poieq}
%\theta^c(t)= c_n\theta([N_t]).
%\end{equation}
%We consider a convergence property of pure step Markov process $\theta^c=\{\theta^c(t);t\ge 0\}$. 
%More precisely, if we denote $M_n^*$ for the law of $\theta^c$ given $x$, 
%we show convergence of $M_n^*$ to a law $M^*$ defined below. 
%Corresponding Ptk $K^n(dh^*;h,x_n)$ in the sense of Section 4b of Chapter XI 
%of \citep{JS} is
%\begin{displaymath}
%r_n\sum_{m=0}^\infty p(dh^*|h,m)p({m}|h,x_n).
%\end{displaymath} 
%\begin{enumerate}
%\item For $m=0,\ldots, n$, and for $\theta= h/c_n$, 
%\begin{displaymath}
%p(m|h,x)=\sum_{y_n\in\{0,1\}^n}1_{\{y^1+\cdots+y^n=m\}}\frac{(1-\theta)f_0(x^j)(1-y^j)+\theta f_\epsilon(x^j)y^j}{(1-\theta)f_0(x^j)+\theta f_\epsilon(x^j)}
%\end{displaymath}
%\item For $m=0,\ldots, n$, $p(A|h,m)=\mathcal{L}(h^*\in A|h,m)$ is defined by 
%$\mathcal{L}( c_n^{-1}(h^*+h)|h,m)=\mathrm{Beta}(\alpha_1+m,\alpha_0+n-m)$. 
%\end{enumerate}
%We will write $n_1=\sum_{i=1}^n y^i$ which is a sufficient statistic for $\theta$. 
%We rewrite the DA procedure, which generate $h_0^n,\ldots$ by 
%\begin{equation}\label{scaledda}
%\left\{\begin{array}{l}h=\theta/c_n,\ \mathbf{simulate}\ y^i|\theta, x^i\sim P(dy|x^i,\theta)\ i=1,\ldots,n),\ 
%\mathbf{set}\ n_1=\sum_{i=1}^n y^i,\\
% \mathbf{simulate}\ \theta^*|x_n,y_n\sim \mathrm{Beta}(\alpha_1+n_1,\alpha_0+n-n_1),\ \mathbf{set}\ h^*=\theta^*/c_n.
% \end{array}\right.
%\end{equation}
Note that by using $t_g$, $\mathbb{P}(y^i=1|\theta,x^i)$
and 
$\mathbb{P}(y^i=0|\theta,x^i)$ become
\begin{equation}\label{flip}
\frac{\theta f_\epsilon(x^i)}{p_\theta^\epsilon(x^i)}
=\frac{\theta+t_{p_\theta^\epsilon}(x^i)}{1+t_{p_\theta^\epsilon}(x^i)},\ 
\frac{(1-\theta) f(x^i)}{p_\theta^\epsilon(x^i)}
=\frac{1-\theta}{1+t_{p_\theta^\epsilon}(x^i)},
\end{equation}
with respectively. 
Throughout in this section, we set $\theta=h/c_n$ and consider uniform convergence property with respect to $h$.  We use the both notation $\theta$ and $h$. Write $X_\theta=o_{\mathbb{P}_n}(\delta_n)$ or $X_\theta=O_{\mathbb{P}_n}(\delta_n)$
for some sequence $\delta_n>0$  
if $\sup_{\theta\in [0,H/c_n]}|X_\theta|$ is 
$o_{\mathbb{P}_n}(\delta_n)$ or $O_{\mathbb{P}_n}(\delta_n)$ 
for any $H>0$ with respectively. 
%Let  
%$\{X_t;t\ge 0\}$ be a stochastic process on a stochastic base $(\Omega,\mathcal{F},\mathbf{F},P_z)$
%with coefficients $(b(\cdot,z),c(\cdot,z),K(\cdot,z)\equiv 0)$ 
%in a sense of Definition 2.18 of Chapter III of \citep{JS}
%and the initial distribution $P^*(\cdot|d=z)$
%defined in (\ref{asymppost}). 
%Set $M^*(\cdot;z)=\mathcal{L}(\{X_t;t\ge 0\}|P_z)$ and $\mathcal{L}(z)=P:=N(0,I)$. 

\begin{lemma}
Let $n_1=\sum_{i=1}^ny^i$ such that 
$y^i|\theta,x^i\sim P(dy|\theta=h/c_n,x^i)$. Under Assumption \ref{Assumption}, 
\begin{eqnarray}\label{eq8}
\mathbb{E}_n[n_1|\theta,x_n]-n\theta&=&h n^{-1/2}\sum_{i=1}^n d(x^i)-h^2I+o_{\mathbb{P}_n}(1), \\
\label{eq24}
r_n^{-1}\mathbb{V}_n[n_1|\theta,x_n]&=&h + o_{\mathbb{P}_n}(1). 
\end{eqnarray}
In particular, both $\mathbb{E}_n[n_1|\theta,x_n]$
and $\mathbb{V}_n[n_1|\theta,x_n]$ are
$O_{\mathbb{P}_n}(r_n)$. 
Moreover, for $k=3,4$, 
\begin{eqnarray}\label{eq7}
\mathbb{E}_n[(n_1-\mathbb{E}_n[n_1|\theta,x_n])^k|\theta,x_n]&=&
O_{\mathbb{P}_n}(r_ n^{k-2}).
\end{eqnarray}
\end{lemma}
\begin{proof}
Note 
\begin{equation}\nonumber
\frac{\theta+x}{1+x}=\theta+(1-\theta)\Big\{x-x^2+\frac{x^3}{1+x}\Big\}. 
\end{equation}
Thus, by (\ref{flip}), 
\begin{displaymath}
\mathbb{E}_n[n_1|\theta,x_n]=\sum_{i=1}^n\frac{\theta+t_{p_\theta^\epsilon}(x^i)}{1+t_{p_\theta^\epsilon}(x^i)}
=n\theta+(1-\theta)\{\sum_{i=1}^n t_{p_\theta^\epsilon}(x^i)-\sum_{i=1}^n t_{p_\theta^\epsilon}^2(x^i)+\sum_{i=1}^n\frac{t_{p_\theta^\epsilon}^3(x^i)}{1+t_{p_\theta^\epsilon}(x^i)}\}. 
\end{displaymath}
We already have a similar expansion in Lemma \ref{lem1} (take $R(x)=x/(1+x)$ here), and so for $\theta=h/c_n$, the inside of the bracket is
\begin{equation}\nonumber
h n^{-1/2}\sum_{i=1}^n d(x^i)-h^2 I+o_{\mathbb{P}_n}(1). 
\end{equation}
Thus we obtain (\ref{eq8}). 
To show (\ref{eq24}), observe
\[
\mathbb{V}_n[n_1|\theta,x_n]=
\sum_{i=1}^n
\frac{\theta+t_{p_\theta^\epsilon}(x^i)}{1+t_{p_\theta^\epsilon}(x^i)}
\frac{1-\theta}{1+t_{p_\theta^\epsilon}(x^i)}
=
(1-\theta)\sum_{i=1}^n
\Big\{\theta+\frac{(1-\theta)t_{p_\theta^\epsilon}(x^i)}{1+t_{p_\theta^\epsilon}(x^i)}\Big\}\Big\{
1-\frac{t_{p_\theta^\epsilon}(x^i)}{1+t_{p_\theta^\epsilon}(x^i)}\Big\}. 
\]
The reading term of the above is $n(1-\theta)\theta$, 
and the remaining terms are dominated by
\begin{equation}\nonumber
|\sum_{i=1}^n\frac{t_{p_\theta^\epsilon}(x^i)}{1+t_{p_\theta^\epsilon}(x^i)}|+
\sum_{i=1}^n\frac{t_{p_\theta^\epsilon}^2(x^i)}{(1+t_{p_\theta^\epsilon}(x^i))^2}. 
\end{equation}
For $\theta=h/c_n$, these are $O_{\mathbb{P}_n}(1)$.
This is clear by expansion $x/(1+x)=x-x^2/(1+x)$ with $\sum_{i=1}^n t_{p_{\theta}^\epsilon}(x^i)=O_{\mathbb{P}_n}(1)$
and $\max_i|t_{p_{\theta}^\epsilon}(x^i)|=o_{\mathbb{P}_n}(1)$. 
This proves (\ref{eq24}). 

Next we show (\ref{eq7}). If $X^i\sim\mathrm{Bernoulli}(p^i)$, 
$|\mathbb{E}(X^i-p^i)^k|=|p^i(1-p^i)^k+(-p^i)^k(1-p^i)|\le 2p^i=2\mathbb{E}[X^i]$. 
Let $S=\sum_{i=1}^n X_i$. Then
\begin{equation}\nonumber
\begin{array}{l}|\mathbb{E}[(S-\mathbb{E}[S])^3]|=
|\sum_{i=1}^n\mathbb{E}[(X^i-\mathbb{E}[X^i])^3]|\le 2\mathbb{E}[S],\\
|\mathbb{E}[(S-\mathbb{E}[S])^4]|=
|\sum_{i=1}^n\mathbb{E}[(X^i-\mathbb{E}[X^i])^4]|+\sum_{i\neq j}\mathbb{V}[X^i]\mathbb{V}[X^j]
\le 2\mathbb{E}[S]+
\mathbb{V}[S]^2\le 
 2\mathbb{E}[S]+
4\mathbb{E}[S]^2.
\end{array}
\end{equation}
Using this fact, 
\begin{equation}\nonumber
\begin{array}{l}
|\mathbb{E}_n[(n_1-\mathbb{E}_n[n_1|\theta,x_n])^3|\theta,x_n]|
\le 2\mathbb{E}_n[n_1|\theta,x_n]=O_{\mathbb{P}_n}(r_n),\\
\mathbb{E}_n[(n_1-\mathbb{E}_n[n_1|\theta,x_n])^4|\theta,x_n]\le 2\mathbb{E}_n[n_1|\theta,x_n]+4\mathbb{E}_n[n_1|\theta,x_n]^2=O_{\mathbb{P}_n}(r_n^2). 
\end{array}
\end{equation}

\end{proof}

\begin{proposition}[Convergence of the coefficients]\label{WCGlem1}
Let $\theta^*$ be the output of (\ref{damix}) when $\theta$ is the previous value. 
Then under Assumption \ref{Assumption}, 
\begin{eqnarray}
n\mathbb{E}_n[\theta^*-\theta|\theta,x_n]&=&\alpha_1+hn^{-1/2}\sum_{i=1}^n d(x^i)-h^2I+o_{\mathbb{P}_n}(1),\label{eq31}\\
nc_n\mathbb{V}_n[\theta^*|\theta,x_n]&=&2h+o_{\mathbb{P}_n}(1)\label{eq32}
\end{eqnarray}
and 
\begin{equation}\label{eq5}
n^2c_n^2\mathbb{E}_n[(\theta^*-\mathbb{E}_n[\theta^*|\theta,x_n])^4|\theta,x_n]=O_{\mathbb{P}_n}(1). 
\end{equation}
\end{proposition}

\begin{proof}
%For notational simplicity, write
%$E(f(h^*)|h,m)=\int_{\mathbf{R}^+} p(dh^*|h,m)f(h^*)$
%and 
%$E(f(n_1)|h,x)=\sum_{m=0}^n p(n_1|h,x)f(n_1)$. 
Since $\theta^*|x_n,y_n$ follows $\mathrm{Beta}(\alpha_1+n_1,\alpha_0+n-n_1)$, we have
\begin{eqnarray*}\nonumber%\label{eq20}
\mathbb{E}_n[\theta^*-\theta|\theta,x_n]&=&
\mathbb{E}_n[\mathbb{E}_n[\theta^*|x_n,y_n]|\theta,x_n]-\theta=\frac{\alpha_1+\mathbb{E}_n[n_1|\theta,x_n]}{\alpha_1+\alpha_0+n}- \theta, \\
&=&\Big(1-\frac{\alpha_1+\alpha_0}{\alpha_1+\alpha_0+n}\Big)\frac{\{\alpha_1+\mathbb{E}_n[n_1|\theta,x_n]-n\theta\}-\theta(\alpha_1+\alpha_0)}{n}
\end{eqnarray*}
and thus (\ref{eq31}) follows by (\ref{eq8}). 
%\begin{displaymath}
%b_n(h,x_n)=\alpha_1+E(n_1|h,x_n)-
%r_nh+o_{P_n}(1)=\alpha_1+hZ_n-h^2I+o_{P_n}(1).
%\end{displaymath}
Observe that 
\begin{equation}\label{eq50}
\mathbb{E}_n[(\theta^*-\theta)^2|x_n,y_n]=
(\mathbb{E}_n[\theta^*|x_n,y_n]-\theta)^2+\mathbb{V}_n[\theta^*|x_n,y_n].
\end{equation}
Expectation of the first term in the right-hand side is
\begin{eqnarray*}
\mathbb{E}_n[(\mathbb{E}_n[\theta^*|x_n,y_n]-\theta)^2|\theta,x_n]
&=&
\mathbb{E}_n[(\mathbb{E}_n[\theta^*|x_n,y_n]-\mathbb{E}_n[\theta^*|\theta,x_n])^2|\theta,x_n]+\mathbb{E}_n[\theta^*-\theta|\theta,x_n]^2\\
&=&
\frac{\mathbb{V}_n[n_1|\theta,x_n]}{(\alpha_1+\alpha_0+n)^2}+\mathbb{E}_n[\theta^*-\theta|\theta,x_n]^2=h/nc_n+O_{\mathbb{P}_n}(1/n^2)
\end{eqnarray*}
where the last equation comes from (\ref{eq24}, \ref{eq31}).  
Expectation of the second term of the right-hand side of (\ref{eq50}) is 
\begin{eqnarray}
\mathbb{E}_n[\mathbb{V}_n[\theta^*|x_n,y_n]|\theta,x_n]&=&\mathbb{E}_n\Big[\frac{(\alpha_1+n_1)(\alpha_0+n-n_1)}{(\alpha_1+\alpha_0+n)^2(\alpha_1+\alpha_0+1+n)}|\theta,x_n\Big]
=\mathbb{E}_n\Big[\frac{nn_1-n_1^2+O(n)}{n^3+O(n^2)}|\theta,x_n\Big]\nonumber\\
&=&\frac{\mathbb{E}_n[n_1|\theta,x_n]}{n^2}+O_{\mathbb{P}_n}\Big(\frac{r_n^2}{n^3}\Big)=\frac{h}{nc_n}+O_{\mathbb{P}_n}\Big(\frac{r_n^2}{n^3}\Big),\label{appendeq1}
\end{eqnarray} 
where in the third equality, we used $\mathbb{E}_n[n_1^2|\theta,x_n]=\mathbb{E}_n[n_1|\theta,x_n]^2+
\mathbb{V}_n[n_1^2|\theta,x_n]=O_{\mathbb{P}_n}(r_n^2)$, and $n^{-2}=o(r^2_n/n^3)$. 
This proves (\ref{eq32}).

Last we show (\ref{eq5}). 
By Jensen's inequality, 
\begin{eqnarray*}
\mathbb{E}_n[(\theta^*-\mathbb{E}_n[\theta^*|\theta,x_n])^4|\theta,x_n]
\le 
%\mathbb{E}_n[(\theta^*-\mathbb{E}_n[\theta^*|x_n,y_n])^4|\theta,x_n]
%=
\mathbb{E}_n[\mathbb{E}_n[(\theta^*-\mathbb{E}_n[\theta^*|x_n,y_n])^4|x_n,y_n]|\theta,x_n]. 
\end{eqnarray*}
Recall that if $X\sim\mathrm{Beta}(\alpha,\beta)$, the kurtosis is
\begin{equation}\nonumber
\frac{\mathbb{E}[(X-\mathbb{E}[X])^4]}{\mathbb{V}[X]^2}
=3+6\frac{(\alpha-\beta)^2(\alpha+\beta+1)-\alpha\beta(\alpha+\beta+2)}{\alpha\beta(\alpha+\beta+2)(\alpha+\beta+3)}
\le 3+\frac{6}{\min\{\alpha,\beta\}}
\end{equation}
where in the second inequality, we used simple fact such as $|\alpha-\beta|\le\max\{\alpha,\beta\}$. 
Since $\theta^*|x_n,y_n\sim\mathrm{Beta}(\alpha_1+n_1,\alpha_0+n-n_1)$ and 
$\min\{\alpha_1+n_1,\alpha_0+n-n_1\}\ge\min\{\alpha_1,\alpha_0\}>0$, we have
\begin{displaymath}
\mathbb{E}_n[(\theta^*-\mathbb{E}_n[\theta^*|x_n,y_n])^4\le 
C\mathbb{V}_n[\theta^*|x_n,y_n]^2
\end{displaymath}
for some constant $C>0$. Hence (\ref{eq5}) follows by (\ref{appendeq1}). 
\end{proof}

Recall that  $\theta^n_0, \theta^n_1,\ldots$ are scaled to $h^n_m=c_n\theta^n_m$. Furthermore, we introduce 
an interpolated process
\begin{equation}\nonumber
X^n_t=h^n_{[r_nt]}=c_n\theta^n_{[r_nt]}
\end{equation}
where $[x]$ is the integer part of $x$. 
Write $\mathcal{L}(\{X_t^n;t\ge 0\}|x_n)$ for the law of $\{X_t^n;t\ge 0\}$ given $x_n$. 
We show that the convergence of $\{X_t^n;t\ge 0\}$ to 
$\{X_t;t\ge 0\}$. There are many studies for the convergence of Markov chain to Markov process. 
See Section 4.8 of \citet{citeulike:5285887} and references therein. We apply Theorem IX.4.21 of \citet{JS} that shows 
the convergence of pure jump Markov process to a diffusion process.  Note that 
$\{X_t^n;t\ge 0\}$ is not a Markov process, since it jumps at deterministic time, but still we can apply the theorem by Proposition VI.6.37 of \citet{JS}. 

\begin{theorem}[Convergence of the DA procedure to a diffusion process]\label{maincor}
By Assumption \ref{Assumption}, $\mathcal{L}(\{X_t^n;t\ge 0\}|x_n)$ 
tends  to $\mathcal{L}(\{X_t;t\ge 0\}|z)$ in distribution, where $z\sim N(0,I)$. 
\end{theorem}

\begin{proof}
By Skorohod's representation theorem, we may assume that there exists a probability space $(\Omega,\mathcal{F},\mathbb{P})$ such that  $z_n(\omega)\rightarrow z(\omega)$ for any $\omega\in\Omega$. 
Let $b(\xi;z)=\alpha_1+\xi z-\xi^2I$ and $c(\xi;z)^2=2\xi$, and set 
\begin{equation}\nonumber
\left\{\begin{array}{l}b^n(\xi;x_n)=\lambda_nc_n\mathbb{E}_n[\theta^*-\theta|\theta=c_n\xi,x_n],\ 
c^n(\xi;x_n)=\lambda_nc_n^2\mathbb{V}_n[\theta^*|\theta=c_n\xi,x_n]. \\
d^n(\xi;x_n)=\lambda_n c_n^4\mathbb{E}_n[(\theta^*-\mathbb{E}_n[\theta^*|\theta=c_n\xi,x_n])^4|\theta=c_n\xi,x_n]. 
\end{array}\right. 
\end{equation}
In Proposition \ref{WCGlem1}, we proved that 
\begin{equation}\nonumber
\mathbb{P}(\sup_{|\xi|<H}|b^n(\xi;x_n)-b(\xi;z_n)|+|c^n(\xi;x_n)-c(\xi;z_n)|+|d^n(\xi;x_n)|>\epsilon)\rightarrow 0
\end{equation}
for any $H>0$, $\epsilon>0$. Thus by local uniform continuity of $b(\xi;z)$ and $c(\xi;z)$ in $(\xi, z)$, we have
$U_n=o_{\mathbb{P}}(1)$ where 
\begin{equation}\nonumber
U_n=\sum_{m=1}^\infty 2^{-m}\min\{1,\sup_{|\xi|<m}|b^n(\xi;x_n)-b(\xi;z)|+|c^n(\xi;x_n)-c(\xi;z)|+|d^n(\xi;x_n)|\}.
\end{equation}
Thus 
by Theorem IX.4.21 of \citet{JS}, 
$\mathcal{L}(\{X_t^n;t\ge 0\}|x_n)$ converges to $\mathcal{L}(\{X_t;t\ge 0\}|z)$ in probability. 
Indeed, by Skorohod's representation theorem again, we may assume $z_n\rightarrow z$ and $U_n\rightarrow 0$
on a probability space $(\Omega',\mathcal{F}',\mathbb{P}')$. 
Then $b^n(\xi;x_n(\omega))\rightarrow b(\xi;z(\omega))$, 
$c^n(\xi;x_n(\omega))\rightarrow c(\xi;z(\omega))$
and 
$d^n(\xi;x_n(\omega))\rightarrow 0$, local uniformly in $\xi$ for any $\omega\in\Omega'$. 
\end{proof}

\bibliographystyle{plainnat}

\end{document}